\documentclass[11pt]{article}

\usepackage{amsmath, amssymb, amsthm}
\usepackage{graphicx}
\usepackage{hyperref}
\usepackage{geometry}
\title{Explicit Formulas for $\mu$-Bases of Planar Rational Quartic Curves}
\author{Weizhen Han\quad\quad  Weikun Sun \\
      {\small School of Science, Tianjin University of Technology and Education, Tianjin 300222, China} \\
        \texttt{0221241093,sunweikun@tute.edu.cn}}
\date{\today}

\newtheorem{theorem}{Theorem}[section]

\newtheorem{lemma}[theorem]{Lemma}

\newtheorem{definition}[theorem]{Definition}
\newtheorem{example}[theorem]{Example}

\newcommand*\Section[2][]{%
  \section{#2}
  \abovedisplayskip=4.0pt plus 2.0pt minus 2.0pt
  \belowdisplayskip=4.0pt plus 2.0pt minus 2.0pt
  \baselineskip 15.3pt
}

\begin{document}

\maketitle

\begin{abstract}
The $\mu$-basis is an algebraic tool originating from the theory of moving curves and moving surfaces, 
and it is widely used in the study of rational curves and surfaces. 
In this paper, we give the explicit formulas for the $\mu$-basis of planar quartic rational parametric curves 
based on redefined vector polynomials, and several illustrative examples are provided. 
Meanwhile, we also discuss the corresponding cases for quadratic and cubic curves.
\end{abstract}

\Section{Introduction}\label{sec:1}
\setcounter{section}{1}\setcounter{equation}{0}
Implicitization of parametric equations of curves and surfaces is an important problem of long history 
in computer graphics and computer-aided geometric design. 
Several tools that are used to implicitize parametric curves and surfaces are resultant, 
Gr\"obner basis and Wu's method. 
In 1995, Sederberg and Chen proposed a new implicitization method, 
the method of moving curves and moving surfaces \cite{sedergerg1995implicitization}. 
On this basis, Cox et al. first introduced the concept of $\mu$-basis, 
which provides a new approach to computing the implicit equation of a planar rational parametric curve \cite{cox1998moving_line}. 
Due to its algebraic and geometric properties, $\mu$-basis has become a bridge connecting the parametric representation 
and the implicit representation of curves and surfaces. By using $\mu$-basis, new methods have been developed 
for studying geometric properties of curves, including point inversion, reparametrization and singularity computation
 \cite{jia2017survey}. Subsequent studies have extended the concept of $\mu$-bases to rational ruled surfaces, 
 tensor-product rational surfaces, and various types of special surfaces \cite{chen2001mu_basis_ruled,chen2003reparametrization,chen2005mu_basis_implicitization,chen2007quadratic_cubic,wang2008quadratic_surfaces,wang2012steiner,shi2012revolution,chen2008singular_points}.

Meanwhile, research into methods for computing $\mu$-bases has continued to advance.
In 2003, Chen and Wang applied vector elimination to the module of moving lines of curves 
and proposed an efficient algorithm for computing $\mu$-bases of plane rational parametric curves \cite{chen2003mu_basis}.
Song and Goldman generalized this method to compute $\mu$-bases of rational parametric curves 
in arbitrary dimensions \cite{song2009mu_bases}.\ Subsequently, Hong et al. proposed another efficient algorithm 
based on the echelon form computation of numerical matrices \cite{hong2017algorithm_mu_bases}.
Huang and Chen extended the notion of $\mu$-bases to univariate polynomial matrices 
and presented an algorithm for computing $\mu$-bases via polynomial matrix factorization \cite{huang2024algorithm}.
Wang observed that these algorithms do not focus on the explicit forms of $\mu$-bases, 
and in \cite{wang2016explicit_mu_bases} gave explicit formulas for the $\mu$-bases of plane rational quadratic and cubic curves.

Motivated by this, with the aid of the vector polynomials $\mathbf{Q}_i$ 
defined by Goldman et al. in \cite{goldman1984vector_elimination}, 
we establish the relationship between $\mathbf{Q}_i \times \mathbf{F}(t)$ and the explicit formulas for $\mu$-bases.
In Section~3, we present explicit formulas for the $\mu$-bases of plane rational quartic curves in the form of theorems.
When the degree is two or three, the explicit formulas for the $\mu$-bases coincide with 
the results of Wang \cite{wang2016explicit_mu_bases}.

\section{Basic knowledge}\label{sec:2}

\subsection{$\mu$-basis}\label{subsec:2.1}
\setcounter{section}{2}\setcounter{equation}{0}
In the real plane $\mathbb{R}^2$, a rational parametric curve is defined by
\[
(x,y)=\left(\frac{f_1(t)}{f_0(t)},\frac{f_2(t)}{f_0(t)}\right),
\]
where $f_0,f_1,f_2$ are polynomials with real coefficients and
$\mathrm{gcd}(f_0,f_1,f_2)=1$.
For convenience, we often consider its homogeneous form in the projective plane $\mathbb{PR}^2$:
\begin{equation}\label{eq:2.1}
	\mathbf{F}(t)=(f_0(t),f_1(t),f_2(t)).
\end{equation}

The degree of $\mathbf{F}(t)$ is defined as
\[
\mathrm{deg}(\mathbf{F}(t))\triangleq
\mathrm{max}\bigl(\mathrm{deg}(f_0),\mathrm{deg}(f_1),\mathrm{deg}(f_2)\bigr).
\]

A moving line $L(x,y,w;t)=0$ in $\mathbb{PR}^2$ is a family of lines corresponding to the parameter $t$:
\begin{equation}\label{eq:2.2}
	L(x,y,w;t)=A(t)x+B(t)y+C(t)w=0,
\end{equation}
where $A(t),B(t),C(t)\in\mathbb{R}[t]$ and
\[
\mathrm{deg}(L(x,y,w;t))=
\mathrm{max}\bigl(\mathrm{deg}(A(t)),\mathrm{deg}(B(t)),\mathrm{deg}(C(t))\bigr).
\]
For convenience, the moving line in \eqref{eq:2.2} is denoted by
$\mathbf{L}(t)=(A(t),B(t),C(t))$.
If
\begin{equation}\label{eq:2.3}
	\mathbf{L}(t)\cdot\mathbf{F}(t)
	=A(t)f_0(t)+B(t)f_1(t)+C(t)f_2(t)\equiv0,
\end{equation}
then $\mathbf{L}(t)$ is said to follow $\mathbf{F}(t)$.
From an algebraic point of view, $\mathbf{L}(t)$ is a syzygy of $\mathbf{F}(t)$.
We denote by $\mathbf{L}$ the set of all moving lines following $\mathbf{F}(t)$.
Then $\mathbf{L}$ is the syzygy module of $\mathbf{F}(t)$.
The $\mu$-basis of a plane rational curve is defined as follows \cite{wang2016explicit_mu_bases}.

\begin{definition}
	
	Let $\mathbf{u}_1(t),\mathbf{u}_2(t)$ be two syzygies of the plane rational curve $\textbf{F}(t)$. If the following two conditions are satisfied:
	
	$1)\;\mathrm{deg}(\mathbf{u}_1(t))+\mathrm{deg}(\mathbf{u}_2(t))=\mathrm{deg}(\mathbf{F}(t))$,
	
	$2)\;\mathbf{u}_1(t)\times\mathbf{u}_2(t)=\lambda\mathbf{F}(t)$, where $\lambda\in \mathbb{R}^*$.\\
	Then $\textbf{u}_1(t),\textbf{u}_2(t)$ are called a $\mu$-basis of $\textbf{F}(t)$.
	
\end{definition}

\begin{lemma}
	
	The implicit equation of the plane rational curve $\mathbf{F}(t)$ can be expressed as:
	\[
	F(x, y, w) = \operatorname{Res}_t\bigl(\mathbf{u}_1(t)\cdot(x, y, w),\; \mathbf{u}_2(t)\cdot(x, y, w)\bigr) = 0,
	\]
	where $\textbf{u}_1(t),\textbf{u}_2(t)$ form a $\mu$-basis of $\textbf{F}(t)$.
	
\end{lemma}

\subsection{Vector Polynomials $\mathbf{Q}_k(t)$}\label{subsec:2.2}

First, we write $\eqref{eq:2.1}$ in the form of a vector polynomial, namely $\mathbf{F}(t)=\displaystyle\sum_{0}^{n}\textbf{p}_i t^i$,
where $\textbf{p}_i\in\mathbb{R}^3$ and $\textbf{p}_0,\textbf{p}_n\ne \textbf{0}.$

For $k=0,1,\cdots,n-1$, $\mathbf{F}(t)$ can be decomposed as follows:
\[
\mathbf{F}(t)=t^{n-k}(\mathbf{p}_nt^k+\cdots+\mathbf{p}_{n-k})+(\mathbf{p}_{n-k-1}t^{n-k-1}+\cdots+\mathbf{p}_0),
\]
Let $\mathbf{Q}_k(t)=\mathbf{p}_nt^k+\cdots+\mathbf{p}_{n-k}$. By the properties of the cross product, we have
\[
\mathbf{Q}_k(t)\times\mathbf{F}(t)=(\mathbf{p}_n\times\mathbf{p}_{n-k-1})t^{n-1}+\text{lower-degree terms},
\]
and $\left(\mathbf{Q}_k\times\mathbf{F}(t)\right)\cdot\mathbf{F}(t)=0$,
which shows that $\mathbf{Q}_0\times\mathbf{F}(t)$，$\mathbf{Q}_1\times\mathbf{F}(t)$，
$\cdots$，$\mathbf{Q}_{n-1}\times\mathbf{F}(t)$are $n$ moving lines of degree $n-1$ following $\mathbf{F}(t)$.

In order to find a $\mu$-basis of the curve, we consider $\mathbf{Q}_0\times\mathbf{F}(t)$ and $\mathbf{Q}_{n-1}\times\mathbf{F}(t)$.

Using the vector cross product identity $(\mathbf{A}\times\mathbf{B})\times(\mathbf{C}\times\mathbf{D})=\mathbf{C}((\mathbf{A}\times\mathbf{B})\cdot\mathbf{D})-\mathbf{D}((\mathbf{A}\times\mathbf{B})\cdot\mathbf{C})$, we have
\begin{equation}\label{eq:2.4}
	\bigl(\mathbf{Q}_{n-1}\times\mathbf{F}(t)\bigr)\times\bigl(\mathbf{Q}_0\times\mathbf{F}(t)\bigr)=\left[\frac{(\mathbf{p}_n\times\mathbf{p}_0)\cdot\mathbf{F}(t)}{t}\right]\mathbf{F}(t).
\end{equation}

\noindent 1) When $\mathbf{p}_n\times\mathbf{p}_0=\mathbf{0}$, $\mathbf{Q}_{n-1}\times\mathbf{F}(t)$ and $\mathbf{Q}_0\times\mathbf{F}(t)$ are linearly dependent over $\mathbb{R}[t]$.\\
\noindent 2) When $\mathbf{p}_n\times\mathbf{p}_0\ne\mathbf{0}$, $\mathbf{Q}_{n-1}\times\mathbf{F}(t)$ and $\mathbf{Q}_0\times\mathbf{F}(t)$ are linearly independent over $\mathbb{R}[t]$, where $\displaystyle\frac{(\mathbf{p}_n\times\mathbf{p}_0)\cdot\mathbf{F}(t)}{t}$ is a scalar polynomial of degree $n-2$, and it can be rewritten as
\[
\frac{(\mathbf{p}_n\times\mathbf{p}_0)\cdot\mathbf{F}(t)}{t}=-\frac{(\mathbf{F}(t)\times\mathbf{p}_0)\cdot\mathbf{p}_n}{t}=-\left(\mathbf{Q}_{n-1}\times\mathbf{F}(t)\right)\cdot \mathbf{p}_n.
\]
If $\mathbf{Q}_{n-1}\times\mathbf{F}(t)$ can be factorized as
\[
\mathbf{Q}_{n-1}\times\mathbf{F}(t)=-\left[\frac{(\mathbf{p}_n\times\mathbf{p}_0)\cdot\mathbf{F}(t)}{t}\right]\cdot\mathbf{g}(t),
\]
where $\mathbf{g}(t)$ is a vector polynomial of degree one satisfying $\mathbf{g}(t)\cdot\mathbf{p}_n=1$,
then $\mathbf{g}(t)\times\left(\mathbf{Q}_0\times\mathbf{F}(t)\right)=\mathbf{F}(t)$ and
$\mathrm{deg}(\mathbf{g}(t))+\mathrm{deg}\left(\mathbf{Q}_0\times\mathbf{F}(t)\right)=\mathrm{deg}(\mathbf{F}(t))$, which shows that $\mathbf{g}(t)$ and $\mathbf{Q}_0\times\mathbf{F}(t)$
form a $\mu$-basis of the curve $\mathbf{F}(t)$.

\section{Explicit Formulas for the $\mu$-Bases of Quartic Curves}\label{sec:3}
\setcounter{section}{3}\setcounter{equation}{0}
For a quartic non-degenerate curve $\mathbf{F}(t)=\mathbf{p}_0+\mathbf{p}_1t+\mathbf{p}_2t^2+\mathbf{p}_3t^3+\mathbf{p}_4t^4$,
where $\mathbf{p}_i\in\mathbb{R}^3,i=0,1,2,3,4$, $\mathbf{p}_0\ne\mathbf{0}$ and $\mathbf{p}_4\ne\mathbf{0}$.
In what follows, $|\mathbf{p}_i\mathbf{p}_j\mathbf{p}_k|$ denotes the determinant formed by the vectors $\mathbf{p}_i^\mathrm{T}$, $\mathbf{p}_j^\mathrm{T}$, and ${\mathbf{p}}_k^\mathrm{T}$.

As stated in Section 2.2, let $\mathbf{Q}_0=\mathbf{p}_4,\mathbf{Q}_1=\mathbf{p}_3+\mathbf{p}_4t,\mathbf{Q}_2=\mathbf{p}_2+\mathbf{p}_3t+\mathbf{p}_4t^2,\mathbf{Q}_3=\mathbf{p}_1+\mathbf{p}_2t+\mathbf{p}_3t^2+\mathbf{p}_4t^3$, then we have:
\begin{eqnarray*}
	\mathbf{Q}_0\times\mathbf{F}(t)&=&(\mathbf{p}_4\times \mathbf{p}_0)+(\mathbf{p}_4\times\mathbf{p}_1)t+(\mathbf{p}_4\times\mathbf{p}_2)t^2+(\mathbf{p}_4\times\mathbf{p}_3)t^3,\\
	\mathbf{Q}_1\times\mathbf{F}(t)&=&(\mathbf{p}_3\times\mathbf{p}_0)+\big[(\mathbf{p}_3\times\mathbf{p}_1)+(\mathbf{p}_4\times\mathbf{p}_0)\big]t+\big[(\mathbf{p}_3\times\mathbf{p}_2)+(\mathbf{p}_4\times\mathbf{p}_1)\big]t^2+(\mathbf{p}_4\times\mathbf{p}_2)t^3,\\
	\mathbf{Q}_2\times\mathbf{F}(t)&=&(\mathbf{p}_2\times\mathbf{p}_0)+\big[(\mathbf{p}_2\times\mathbf{p}_1)+(\mathbf{p}_3\times\mathbf{p}_0)\big]t+\big[(\mathbf{p}_3\times\mathbf{p}_1)+(\mathbf{p}_4\times\mathbf{p}_0)\big]t^2+(\mathbf{p}_4\times\mathbf{p}_1)t^3,\\
	\mathbf{Q}_3\times\mathbf{F}(t)&=&(\mathbf{p}_1\times\mathbf{p}_0)+(\mathbf{p}_2\times\mathbf{p}_0)t+(\mathbf{p}_3\times \mathbf{p}_0)t^2+(\mathbf{p}_4\times\mathbf{p}_0)t^3.
\end{eqnarray*}
In the following, we classify the discussion according to whether $\mathbf{p}_4\times\mathbf{p}_0$ is the zero vector.

\subsection{The Case $\mathbf{p}_4\times\mathbf{p}_0=\textbf{0}$}\label{subsec:3.1}
Since in this case $\mathbf{Q}_3\times\mathbf{F}(t)$ and $\mathbf{Q}_0\times\mathbf{F}(t)$ are linearly dependent over $\mathbb{R}[t]$,
we consider the simplified $\displaystyle\frac{\mathbf{Q}_0\times\mathbf{F}(t)}{t}$ and $\mathbf{Q}_1\times\mathbf{F}(t)$,
and have
$$\displaystyle\frac{\mathbf{Q}_0\times\mathbf{F}(t)}{t}=(\mathbf{p}_4\times\mathbf{p}_1)+(\mathbf{p}_4\times\mathbf{p}_2)t+(\mathbf{p}_4\times\mathbf{p}_3)t^2,$$
$$\displaystyle\frac{\mathbf{Q}_0\times\mathbf{F}(t)}{t}\times(\mathbf{Q}_1\times\mathbf{F}(t))=-(|\mathbf{p}_1\mathbf{p}_3\mathbf{p}_4|+|\mathbf{p}_2\mathbf{p}_3\mathbf{p}_4|t)\cdot\mathbf{F}(t).$$
\noindent
and when $|\mathbf{p}_1\mathbf{p}_3\mathbf{p}_4|^2-|\mathbf{p}_1\mathbf{p}_2\mathbf{p}_4|\cdot|\mathbf{p}_2\mathbf{p}_3\mathbf{p}_4|=0$,
$\displaystyle\frac{\mathbf{Q}_0\times\mathbf{F}(t)}{t}$ can be decomposed as follows
\[
\displaystyle\frac{\mathbf{Q}_0\times\mathbf{F}(t)}{t}=\left(|\mathbf{p}_1\mathbf{p}_3\mathbf{p}_4|+|\mathbf{p}_2\mathbf{p}_3\mathbf{p}_4|t\right)\left(\displaystyle\frac{\mathbf{p}_4\times\mathbf{p}_1}{|\mathbf{p}_1\mathbf{p}_3\mathbf{p}_4|}+\displaystyle\frac{\mathbf{p}_4\times\mathbf{p}_3}{|\mathbf{p}_2\mathbf{p}_3\mathbf{p}_4|}t\right).
\]
Therefore, let $\mathbf{L}(t)=\displaystyle\frac{\mathbf{p}_4\times\mathbf{p}_1}{|\mathbf{p}_1\mathbf{p}_3\mathbf{p}_4|}+\displaystyle\frac{\mathbf{p}_4\times\mathbf{p}_3}{|\mathbf{p}_2\mathbf{p}_3\mathbf{p}_4|}t$,
then we have $\mathbf{L}(t)\times(\mathbf{Q}_1\times\mathbf{F}(t))=-\mathbf{F}(t)$, and the degrees of $\mathbf{L}(t)$ and $\mathbf{Q}_1\times\mathbf{F}(t)$ are 1 and 3, respectively.

When $|\mathbf{p}_1\mathbf{p}_3\mathbf{p}_4|^2-|\mathbf{p}_1\mathbf{p}_2\mathbf{p}_4|\cdot|\mathbf{p}_2\mathbf{p}_3\mathbf{p}_4|\ne0$,
consider the vector identity in $\mathbb{R}^3$
\[
|\mathbf{p}_2\mathbf{p}_3\mathbf{p}_4|\cdot\mathbf{p}_1-|\mathbf{p}_1\mathbf{p}_3\mathbf{p}_4|\cdot\mathbf{p}_2+|\mathbf{p}_1\mathbf{p}_2\mathbf{p}_4|\cdot\mathbf{p}_3-|\mathbf{p}_1\mathbf{p}_2\mathbf{p}_3|\cdot\mathbf{p}_4=\mathbf{0},
\]
taking the cross product with $\mathbf{p}_4$ on both sides, we obtain
\begin{equation}\label{eq:3.1}
	|\mathbf{p}_2\mathbf{p}_3\mathbf{p}_4|\cdot(\mathbf{p}_4\times\mathbf{p}_1)-|\mathbf{p}_1\mathbf{p}_3\mathbf{p}_4|\cdot(\mathbf{p}_4\times\mathbf{p}_2)+|\mathbf{p}_1\mathbf{p}_2\mathbf{p}_4|\cdot(\mathbf{p}_4\times\mathbf{p}_3)=\mathbf{0}.
\end{equation}
Using the above formula, we can construct a polynomial of degree 2
\begin{equation}\label{eq:3.2}
	\mathbf{L}(t)=|\mathbf{p}_2\mathbf{p}_3\mathbf{p}_4|\cdot(\mathbf{Q}_2\times\mathbf{F}(t))-|\mathbf{p}_1\mathbf{p}_3\mathbf{p}_4|\cdot(\mathbf{Q}_1\times\mathbf{F}(t))+|\mathbf{p}_1\mathbf{p}_2\mathbf{p}_4|\cdot(\mathbf{Q}_0\times\mathbf{F}(t))
\end{equation}
\noindent It is easy to see that $\mathbf{L}(t)\cdot\mathbf{F}(t)\equiv0$, and
\begin{eqnarray*}
	\displaystyle\frac{\mathbf{Q}_0\times\mathbf{F}(t)}{t}\times\mathbf{L}(t)&=&\displaystyle\frac{\mathbf{L}(t)\cdot\mathbf{p}_4}{t}\mathbf{F}(t)\\
	&=&(|\mathbf{p}_1\mathbf{p}_3\mathbf{p}_4|^2-|\mathbf{p}_1\mathbf{p}_2\mathbf{p}_4|\cdot|\mathbf{p}_2\mathbf{p}_3\mathbf{p}_4|)\mathbf{F}(t),
\end{eqnarray*}
At this point, $\displaystyle\frac{\mathbf{Q}_0\times\mathbf{F}(t)}{t}$ and $\mathbf{L}(t)$ form a $\mu$-basis of the curve $\mathbf{F}(t)$.
In summary, we have
\begin{theorem}\label{theorem3.1}

	If the quartic non-degenerate curve $\mathbf{F}(t)$ satisfies $\mathbf{p}_4\times\mathbf{p}_0=\textbf{0}$，then we have:
	
	$1)$When $|\mathbf{p}_1\mathbf{p}_3\mathbf{p}_4|^2-|\mathbf{p}_1\mathbf{p}_2\mathbf{p}_4|\cdot|\mathbf{p}_2\mathbf{p}_3\mathbf{p}_4|=0$,
	
	$\mathbf{u}(t)=\mathbf{Q}_1\times\mathbf{F}(t)$ and $\mathbf{L}(t)=\displaystyle\frac{\mathbf{p}_4\times\mathbf{p}_1}{|\mathbf{p}_1\mathbf{p}_3\mathbf{p}_4|}+\displaystyle\frac{\mathbf{p}_4\times\mathbf{p}_3}{|\mathbf{p}_2\mathbf{p}_3\mathbf{p}_4|}t$ form a $\mu$-basis of the curve $\mathbf{F}(t)$.
	
	$2)$When $|\mathbf{p}_1\mathbf{p}_3\mathbf{p}_4|^2-|\mathbf{p}_1\mathbf{p}_2\mathbf{p}_4|\cdot|\mathbf{p}_2\mathbf{p}_3\mathbf{p}_4|\ne0$,
	
	$\mathbf{u}(t)=\displaystyle\frac{\mathbf{Q}_0\times\mathbf{F}(t)}{t}$和$\mathbf{L}(t)=\square+\triangle t+\star t^2$ form a $\mu$-basis of the curve $\mathbf{F}(t)$，where
	\begin{eqnarray*}
		\square&=&|\mathbf{p}_2\mathbf{p}_3\mathbf{p}_4|\left(\mathbf{p}_2\times\mathbf{p}_0\right)-|\mathbf{p}_1\mathbf{p}_3\mathbf{p}_4|\left(\mathbf{p}_3\times\mathbf{p}_0\right),\\
		\triangle&=&||\mathbf{p}_2\mathbf{p}_3\mathbf{p}_4|\left[\left(\mathbf{p}_2\times\mathbf{p}_1\right)+\left(\mathbf{p}_3\times\mathbf{p}_0\right)\right]-|\mathbf{p}_1\mathbf{p}_3\mathbf{p}_4|\left(\mathbf{p}_3\times\mathbf{p}_1\right)+|\mathbf{p}_1\mathbf{p}_2\mathbf{p}_4|\left(\mathbf{p}_4\times\mathbf{p}_1\right),\\
		\star&=&|\mathbf{p}_2\mathbf{p}_3\mathbf{p}_4|\left(\mathbf{p}_3\times\mathbf{p}_1\right)-|\mathbf{p}_1\mathbf{p}_3\mathbf{p}_4|\left[\left(\mathbf{p}_3\times\mathbf{p}_2\right)+\left(\mathbf{p}_4\times\mathbf{p}_1\right)\right]+|\mathbf{p}_1\mathbf{p}_2\mathbf{p}_4|\left(\mathbf{p}_4\times\mathbf{p}_2\right).
	\end{eqnarray*}
	
\end{theorem}

\begin{example}
	
	Given a quartic rational parametric curve in the plane, its vector polynomial form is
	$\mathbf{F}(t)=\mathbf{p}_0+\mathbf{p}_1t+\mathbf{p}_2t^2+\mathbf{p}_3t^3+\mathbf{p}_4t^4$，where
	$\mathbf{p}_0 = (1, 0, 0), \mathbf{p}_1 = (0, 1, 0), \mathbf{p}_2 = (0, 0, 1), \mathbf{p}_3 = (1, 1, 1), \mathbf{p}_4 = (2,0,0)$.
	
\end{example}

\indent {Solution}\ \
This curve satisfies $\mathbf{p}_4 \times\mathbf{p}_0 =  \textbf{0}$,
$|\mathbf{p}_1\mathbf{p}_3\mathbf{p}_4|^2-|\mathbf{p}_1\mathbf{p}_2\mathbf{p}_4|\cdot|\mathbf{p}_2\mathbf{p}_3\mathbf{p}_4|=8\ne0$,
and further computation gives
\begin{eqnarray*}
	\square&=&|\mathbf{p}_2\mathbf{p}_3\mathbf{p}_4|\left(\mathbf{p}_2\times\mathbf{p}_0\right)-|\mathbf{p}_1\mathbf{p}_3\mathbf{p}_4|\left(\mathbf{p}_3\times\mathbf{p}_0\right)=(0,-4,2),\\
	\triangle&=&||\mathbf{p}_2\mathbf{p}_3\mathbf{p}_4|\left[\left(\mathbf{p}_2\times\mathbf{p}_1\right)+\left(\mathbf{p}_3\times\mathbf{p}_0\right)\right]-|\mathbf{p}_1\mathbf{p}_3\mathbf{p}_4|\left(\mathbf{p}_3\times\mathbf{p}_1\right)+|\mathbf{p}_1\mathbf{p}_2\mathbf{p}_4|\left(\mathbf{p}_4\times\mathbf{p}_1\right)=(4,-2,4),\\
	\star&=&|\mathbf{p}_2\mathbf{p}_3\mathbf{p}_4|\left(\mathbf{p}_3\times\mathbf{p}_1\right)-|\mathbf{p}_1\mathbf{p}_3\mathbf{p}_4|\left[\left(\mathbf{p}_3\times\mathbf{p}_2\right)+\left(\mathbf{p}_4\times\mathbf{p}_1\right)\right]+|\mathbf{p}_1\mathbf{p}_2\mathbf{p}_4|\left(\mathbf{p}_4\times\mathbf{p}_2\right)=(0,-2,-6).
\end{eqnarray*}
Therefore, we can directly obtain
$$\displaystyle\mathbf{u}(t) =\displaystyle\frac{\mathbf{Q}_0 \times \mathbf{F}(t)}{t} = (0,-2t^2-2t,2t^2+2),$$
$$\mathbf{L}(t)=\square+\triangle t+\star t^2=(4t,-2t^2-2t-4,-6t^2+4t+2).$$
Verification gives $\mathbf{u}(t) \times \mathbf{L}(t) = 8F(t)$, so $\mathbf{u}(t)$ and $\mathbf{L}(t)$ form a $\mu$-basis of the curve.

\noindent Furthermore, the implicit equation of the curve $\mathbf{F}(t)$ is
\begin{align*}
	F(x,y,w) &= \det\begin{pmatrix}
		-4y+2w & 0 & 2w & 0 \\
		4x-2y+4w & -4y+2w & -2y & 2w\\
		-2y-6w & 4x-2y+4w & -2y+2w &-2y\\
		0 & -2y-6w &  0 & -2y+2w
	\end{pmatrix} \\
	&= \begin{aligned}[t]
		&320w^4+128w^3x-448w^3y+64w^2x^2-256w^2xy+256w^2y^2 \\
		&-64wx^2y-64wxy^2+64wy^3+64xy^3+64y^4=0.
	\end{aligned}
\end{align*}

\subsection{The Case $\mathbf{p}_4\times\mathbf{p}_0\ne\textbf{0}$}\label{subsec:3.2}
Similarly to formula $\eqref{eq:3.1}$, we can obtain:
\begin{equation}\label{eq:3.3}
	|\mathbf{p}_b\mathbf{p}_c\mathbf{p}_4|\cdot(\mathbf{p}_4\times\mathbf{p}_a)-|\mathbf{p}_a\mathbf{p}_c\mathbf{p}_4|\cdot(\mathbf{p}_4\times\mathbf{p}_b)+|\mathbf{p}_a\mathbf{p}_b\mathbf{p}_4|\cdot(\mathbf{p}_4\times\mathbf{p}_c)=\mathbf{0}.
\end{equation}

\noindent Using the above formula, we obtain four quadratic moving lines following $\mathbf{F}(t)$,
\begin{eqnarray*}
	\mathbf{u}_0&=&|\mathbf{p}_1\mathbf{p}_2\mathbf{p}_4|(\mathbf{Q}_0\times\mathbf{F}(t))-|\mathbf{p}_1\mathbf{p}_3\mathbf{p}_4|(\mathbf{Q}_1\times\mathbf{F}(t))+|\mathbf{p}_2\mathbf{p}_3\mathbf{p}_4|(\mathbf{Q}_2\times\mathbf{F}(t)),\\
	\mathbf{u}_1&=&|\mathbf{p}_0\mathbf{p}_2\mathbf{p}_4|(\mathbf{Q}_0\times\mathbf{F}(t))-|\mathbf{p}_0\mathbf{p}_3\mathbf{p}_4|(\mathbf{Q}_1\times\mathbf{F}(t))+|\mathbf{p}_2\mathbf{p}_3\mathbf{p}_4|(\mathbf{Q}_3\times\mathbf{F}(t)),\\
	\mathbf{u}_2&=&|\mathbf{p}_0\mathbf{p}_1\mathbf{p}_4|(\mathbf{Q}_0\times\mathbf{F}(t))-|\mathbf{p}_0\mathbf{p}_3\mathbf{p}_4|(\mathbf{Q}_2\times\mathbf{F}(t))+|\mathbf{p}_1\mathbf{p}_3\mathbf{p}_4|(\mathbf{Q}_3\times\mathbf{F}(t)),\\
	\mathbf{u}_3&=&|\mathbf{p}_0\mathbf{p}_1\mathbf{p}_4|(\mathbf{Q}_1\times\mathbf{F}(t))-|\mathbf{p}_0\mathbf{p}_2\mathbf{p}_4|(\mathbf{Q}_2\times\mathbf{F}(t))+|\mathbf{p}_1\mathbf{p}_2\mathbf{p}_4|(\mathbf{Q}_3\times\mathbf{F}(t)).
\end{eqnarray*}
Taking the cross product of any $\mathbf{u}_i,\mathbf{u}_j,(i,j=0,1,2,3)$, and simplifying by using the Grassmann--Pl\"ucker relations, we obtain:
\begin{equation}\label{eq:3.4}
	\mathbf{u}_i\times\mathbf{u}_j=|\mathbf{p}_a\mathbf{p}_b\mathbf{p}_4|\cdot E\cdot\mathbf{F}(t),
\end{equation}
where $0\leq a<b\leq3$ is a pair of indices distinct from $i$ and $j$, and
\begin{eqnarray*}
	E&=&|\mathbf{p}_0\mathbf{p}_2\mathbf{p}_4|^2-(|\mathbf{p}_0\mathbf{p}_2\mathbf{p}_3|+|\mathbf{p}_0\mathbf{p}_1\mathbf{p}_4|)\cdot|\mathbf{p}_0\mathbf{p}_3\mathbf{p}_4|\\
	&&+\,|\mathbf{p}_0\mathbf{p}_1\mathbf{p}_3|\cdot|\mathbf{p}_1\mathbf{p}_3\mathbf{p}_4|-|\mathbf{p}_0\mathbf{p}_1\mathbf{p}_4|\cdot|\mathbf{p}_1\mathbf{p}_2\mathbf{p}_4|-|\mathbf{p}_0\mathbf{p}_1\mathbf{p}_2|\cdot|\mathbf{p}_2\mathbf{p}_3\mathbf{p}_4|.	
\end{eqnarray*}

\noindent
The non-degeneracy of the curve guarantees that there exists some $|\mathbf{p}_a\mathbf{p}_b\mathbf{p}_4|\ne0$,
and hence, when $E\ne0$, according to the definition of a $\mu$-basis, $\mathbf{u}_i$ and $\mathbf{u}_j$ form a $\mu$-basis of the curve $\mathbf{F}(t)$.

In order to consider the case when $E=0$, we first give the following lemma:

	\begin{lemma}\label{lemma3.3}
		Given a non-degenerate quartic curve $\mathbf{F}(t)=\displaystyle\sum_{i=0}^{4}\mathbf{p}it^i$ in $\mathbb{R}^3$,
		let $S$ be the matrix of a non-degenerate linear transformation,
		and let $\mathbf{F'}(t)=\mathbf{F}(t)S=\displaystyle\sum{i=0}^{4}(\mathbf{p}_iS)t^i$ be the transformed curve.
		Let $\mathbf{u'}_1(t)$ and $\mathbf{u'}_2(t)$ be a $\mu$-basis of it.
		Let $\mathbf{u}_1(t)=\mathbf{u'}_1(t)S^{-T}, \mathbf{u}_2(t)=\mathbf{u'}_2(t)S^{-T}$. Then
		$\mathbf{u}_1(t), \mathbf{u}_2(t)$ form a $\mu$-basis of the curve $\mathbf{F}(t)$.
	\end{lemma}

\indent {Proof}\ \  
First, since
\vspace{-0.3cm}
\begin{eqnarray*}
	\mathbf{u}_i(t)\cdot\mathbf{F}(t)&=&\mathbf{u}_i(t)\mathbf{F}(t)^T=\mathbf{u'}_i(t)S^{-T}\mathbf{F}(t)^T=\mathbf{u'}_i(t)(\mathbf{F}(t)S^{-1})^T\\
	&=&\mathbf{u'}_i(t)\mathbf{F'}(t)^T=\mathbf{u'}_i(t)\cdot\mathbf{F'}(t)=0,
\end{eqnarray*}
so we have that $\mathbf{u}_1(t)$ and $\mathbf{u}_2(t)$ follow $\mathbf{F}(t)$.

Since a linear transformation does not change the degree of a polynomial, it remains only to prove that
$\mathbf{u}_1(t)\times\mathbf{u}_2(t)=\lambda\mathbf{F}(t)$.

Consider the cross product extended to the polynomial vector space $\mathbb R[t]^3$.
For any $\mathbf h_1(t),\mathbf h_2(t)\in\mathbb R[t]^3$ and any
invertible matrix $B\in GL_3(\mathbb R)$, we have
\begin{equation}\label{eq:3.5}
	(\mathbf h_1(t)B)\times(\mathbf h_2(t)B)
	=
	\det(B)(\mathbf h_1(t)\times\mathbf h_2(t))B^{-T}.
\end{equation}

Using $\eqref{eq:3.5}$ and $\mathbf{u'}_1(t)\times\mathbf{u'}_2(t)=\lambda'\mathbf{F'}(t)$，$\lambda'\in\mathbb{R}^*$，we obtain
\begin{eqnarray*}
	\mathbf{u}_1(t)\times\mathbf{u}_2(t)&=&(\mathbf{u'}_1(t)S^{-T})\times(\mathbf{u'}_2(t)S^{-T})=\det(S^{-T})(\mathbf{u'}_1(t)\times\mathbf{u'}_2(t))(S^{-T})^{-T}\\
	&=&\frac{1}{\det(S)}\lambda'\mathbf{F'}(t)S=\frac{1}{\det(S)}\lambda'\mathbf{F}(t)S^{-1}S=\lambda\mathbf{F}(t),\ \text{where}\ \lambda=\frac{1}{\det(S)}\lambda'.
\end{eqnarray*}
This completes the proof.

Thus, when $E=0$, without loss of generality, assume that $|\mathbf{p}_0\mathbf{p}_1\mathbf{p}_4|\ne0$. Let the invertible matrix
$S=(\mathbf p_0^T, \mathbf p_1^T, \mathbf p_4^T)^T$ and $\mathbf{F'}(t)=\mathbf{F}(t)S^{-1}$. Then
among $\mathbf{p'}_i=\mathbf{p}_iS^{-1}(i=0,1,2,3,4)$,
$\mathbf{p'}_0$, $\mathbf{p'}_1$, and $\mathbf{p'}_4$ form a standard basis of $\mathbb{R}^3$.

Suppose $\mathbf{p'}_2=(a,b,c)$、$\mathbf{p'}_3=(d,e,f)$，then
$$
\mathbf{F'}(t)=(1+at^2+dt^3,t+bt^2+et^3,ct^2+ft^3+t^4).
$$
Substituting $\mathbf{p'}_i$ into $E$, we obtain.
$$
E=b^2-(bf-ce+1)e-df+a-c(ae-bd).
$$
Let the degree-one moving line following $\mathbf{F'}(t)$ be $\mathbf{L}(t)=\mathbf{A}+\mathbf{B}t=(x_0,y_0,z_0)+(x_1,y_1,z_1)t$, then we have
\[
\mathbf{L}(t)\cdot\mathbf{F'}(t)\equiv0
\iff 
\left\{
\begin{array}{l}
	\left.
	\begin{array}{l}
		x_0 = 0 \\
		z_1 = 0 \\
		y_0 = -x_1 \\
		y_1 = b x_1 - c z_0
	\end{array}
	\right\}
	\implies 
	\mathbf{L}(t) = \left( x_1 t,\ -x_1 + (b x_1 - c z_0) t,\ z_0 \right) \\[15pt]
	\begin{bmatrix}
		a + b^2 - e & f - bc \\
		d + eb & 1 - ce
	\end{bmatrix}
	\begin{bmatrix}
		x_1 \\
		z_0
	\end{bmatrix}
	=
	\begin{bmatrix}
		0 \\
		0
	\end{bmatrix}
\end{array}
\right.
\]
then $\det \begin{bmatrix}
	a + b^2 - e & f - bc \\
	d + eb & 1 - ce
\end{bmatrix}=E$，show:
\[
\small
E=0
\Longleftrightarrow (x_1,z_0)\ne(0,0)
\Longleftrightarrow \exists\ \mathbf{L}(t)\text{ of degree one following }\mathbf{F'}(t)
\Longleftrightarrow Syz(\mathbf{F'})=R(-1)\oplus R(-3).
\]

Using $\mathbf{L}(t)\cdot\mathbf{F'}(t)=(\mathbf{A}+\mathbf{B}t)\cdot(\mathbf{p}_0+\mathbf{p}_1t+\mathbf{p}_2t^2+\mathbf{p}_3t^3+\mathbf{p}_4t^4)=0$，it can be simplified
\begin{eqnarray*}
	(\mathbf{Q}_0\times\mathbf{F'}(t))\times\mathbf{L}(t)&=&(\mathbf{p}_{4}\cdot\mathbf{A})\cdot\mathbf{F'}(t)=z_0\cdot\mathbf{F'}(t),\\
	(\mathbf{Q}_1\times\mathbf{F'}(t))\times\mathbf{L}(t)&=&(\mathbf{p}_{3}\cdot\mathbf{A})\cdot\mathbf{F'}(t)=(fz_0-ex_1)\cdot\mathbf{F'}(t),\\
	(\mathbf{Q}_2\times\mathbf{F'}(t))\times\mathbf{L}(t)&=&(\mathbf{p}_{2}\cdot\mathbf{A})\cdot\mathbf{F'}(t)=(cz_0-bx_1)\cdot\mathbf{F'}(t),\\
	(\mathbf{Q}_3\times\mathbf{F'}(t))\times\mathbf{L}(t)&=&(\mathbf{p}_{1}\cdot\mathbf{A})\cdot\mathbf{F'}(t)=-x_1\cdot\mathbf{F'}(t).
\end{eqnarray*}
Since $(x_1,z_0)\ne(0,0)$, if $(\mathbf{Q}_0\times\mathbf{F'}(t))\times\mathbf{L}(t)=\mathbf{0}$, then $(\mathbf{Q}_3\times\mathbf{F'}(t))\times\mathbf{L}(t)\ne\mathbf{0}$, and vice versa.
Thus, when $E=0$, the curve $\mathbf{F}^{'}(t)$ has a $\mu$-basis whose degrees are 1 and 3, respectively. According to Lemma \ref{lemma3.3}, the original curve $\mathbf{F}(t)$ also has such a $\mu$-basis. In summary, we have

\begin{theorem}\label{theorem3.4}
	If the quartic non-degenerate curve $\mathbf{F}(t)$ satisfies $\mathbf{p}_4\times\mathbf{p}_0\ne\textbf{0}$, then we have:
	
	1) When $E\ne0$, $\mathbf{u}_i$ and $\mathbf{u}_j$ in $\eqref{eq:3.4}$ form a $\mu$-basis of the curve $\mathbf{F}(t)$.
	
	2) When $E=0$, the curve $\mathbf{F}(t)$ has a $\mu$-basis whose degrees are 1 and 3, respectively.
	
\end{theorem}

\indent {Remark}\ \  When $E=0$, all nonzero $\mathbf{u}_i(t)$ are pairwise parallel in $\mathbb{R}^3[t]$. Take
$$
\mathbf{L}(t)=(x_0+x_1t,y_0+y_1t,z_0+z_1t)=\gcd(\mathbf{u}_i)
$$
as a degree-one moving line. If $z_0=0$, then the cubic moving line is $\mathbf{Q}_3\times\mathbf{F}(t)$; if $x_1=0$, then the cubic moving line is $\mathbf{Q}_0\times\mathbf{F}(t)$.

Therefore, when a non-degenerate curve is given, we can obtain a $\mu$-basis of the curve according to the following procedure.
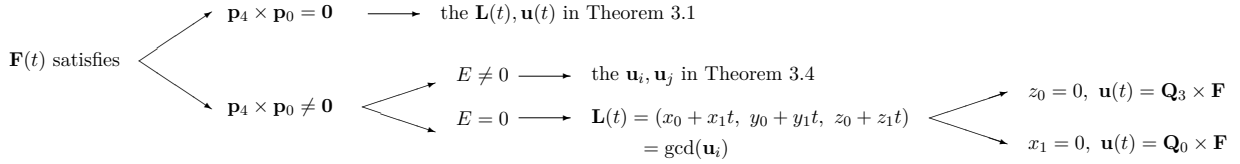
\begin{figure}[htbp]
	\centering
	\setlength{\unitlength}{1mm}
	\resizebox{\textwidth}{!}{%
		\begin{picture}(240,32)
			
			\put(0,17){\makebox(0,0)[l]{$\mathbf F(t)\ \text{satisfies}$}}
			
			\put(25,17){\vector(3,2){14}}
			\put(25,17){\vector(3,-2){14}}
			
			\put(42,26){\makebox(0,0)[l]{$\mathbf p_4\times\mathbf p_0=\mathbf 0$}}
			\put(42,8){\makebox(0,0)[l]{$\mathbf p_4\times\mathbf p_0\neq\mathbf 0$}}
			
			\put(69,26){\vector(1,0){10}}
			\put(83,26){\makebox(0,0)[l]{$\text{the } \mathbf L(t),\mathbf u(t)\ \text{in Theorem \ref{theorem3.1}}$}}
			
			\put(68,8){\vector(3,1){14}}
			\put(68,8){\vector(3,-1){14}}
			
			\put(86,14){\makebox(0,0)[l]{$E\neq 0$}}
			\put(86,6){\makebox(0,0)[l]{$E=0$}}
			
			\put(98,14){\vector(1,0){10}}
			\put(112,14){\makebox(0,0)[l]{$\text{the } \mathbf u_i,\mathbf u_j\ \text{in Theorem \ref{theorem3.4} }$}}
			
			\put(98,6){\vector(1,0){10}}
			\put(112,6){\makebox(0,0)[l]{$\mathbf L(t)=(x_0+x_1t,\ y_0+y_1t,\ z_0+z_1t)$}}
			\put(122,0){\makebox(0,0)[l]{$=\gcd(\mathbf u_i)$}}
			
			\put(177,6){\vector(3,1){15}}
			\put(177,6){\vector(3,-1){15}}
			
			\put(196,11){\makebox(0,0)[l]{$z_0=0,\ \mathbf u(t)=\mathbf Q_3\times\mathbf F$}}
			\put(196,1){\makebox(0,0)[l]{$x_1=0,\ \mathbf u(t)=\mathbf Q_0\times\mathbf F$}}
			
		\end{picture}%
	}
\end{figure}

\begin{example}
	
	Given a plane quartic rational parametric curve $\mathbf{F}(t)=\mathbf{p}_0+\mathbf{p}_1t+\mathbf{p}_2t^2+\mathbf{p}_3t^3+\mathbf{p}_4t^4$，where
	\[
	\mathbf{p}_0 = (1, 0, 0),\quad\mathbf{p}_1 = (1, 0, 0),\quad\mathbf{p}_2 = (1,1,0),\quad \mathbf{p}_3 = (0, 1, 1),\quad\mathbf{p}_4 = (0,0,1).
	\]
	
\end{example}

\indent {Solution}\ \ 
This curve satisfies $\mathbf{p}_4 \times\mathbf{p}_0 \ne \textbf{0}$ and $E=0$. Computing $\mathbf{u}_i$ gives:
\begin{eqnarray*}
	\mathbf{u}_0&=&(0,t^2,-t)=t(0,t,-1),\\
	\mathbf{u}_1&=&(0,-t,1)=-(0,t,-1),\\
	\mathbf{u}_2&=&(0,-t^2-t,t+1)=-(t+1)(0,t,-1),\\
	\mathbf{u}_3&=&(0,-t^2-t,t+1)=-(t+1)(0,t,-1).
\end{eqnarray*}
Let $\mathbf{L}(t)=(0,t,-1)$. Computing $\mathbf{Q}_i\times\mathbf{F}(t)$ gives:
\begin{eqnarray*}
	\mathbf{Q}_0\times\mathbf{F}(t)&=&(-t^2-t^3,1+t+t^2,0),\\
	\mathbf{Q}_1\times\mathbf{F}(t)&=&(-t^2-t^3,1+2t+2t^2+t^3,-1-t-t^2),\\
	\mathbf{Q}_2\times\mathbf{F}(t)&=&(0,t+2t^2+t^3,-1-2t-t^2)=(1+t)^2\mathbf{L}(t),\\
	\mathbf{Q}_3\times\mathbf{F}(t)&=&(0,t^2+t^3,-t-t^2)=t(1+t)\mathbf{L}(t).
\end{eqnarray*}
and
\begin{eqnarray*}
	(\mathbf{Q}_0\times\mathbf{F}(t))\times\mathbf{L}(t)&=&-\mathbf{F}(t),\\
	(\mathbf{Q}_1\times\mathbf{F}(t))\times\mathbf{L}(t)&=&-\mathbf{F}(t).
\end{eqnarray*}
Thus, according to the definition of a $\mu$-basis, both $\mathbf{L}(t)$ and $\mathbf{Q}_0\times\mathbf{F}(t)$, and $\mathbf{L}(t)$ and $\mathbf{Q}_1\times\mathbf{F}(t)$, form a $\mu$-basis of the curve $\mathbf{F}(t)$.
Using $\mathbf{Q}_0\times\mathbf{F}(t)$ and $\mathbf{L}(t)$, the implicit equation of the curve $\mathbf{F}(t)$ is
\begin{align*}
	F(x,y,w) &= \det\begin{pmatrix}
		y & -w & 0 & 0 \\
		0 & y & -w & 0 \\
		0 & 0 & y & -w \\
		-x & y-x & y & y
	\end{pmatrix} \\
	&= \begin{aligned}[t]
		&y^{4} + y^{3}w + y(y-x)w^2 -xw^3= 0.
	\end{aligned}
\end{align*}
	
\section{Results for Quadratic and Cubic Curves}\label{sec:4}

\subsection{Quadratic Non-degenerate Curves}\label{subsec:4.1}

Let $\mathbf{F}(t)=\mathbf{p}_0+\mathbf{p}_1t+\mathbf{p}_2t^2$ be a curve, where $\mathbf{p}_i\in\mathbb{R}^3,i=0,1,2,3$, $\mathbf{p}_0\ne0$ and $\mathbf{p}_2\ne0$.
Let $M=(\mathbf{p}_0^T,\mathbf{p}_1^T,\mathbf{p}_2^T)$ be a $3\times3$ matrix. Then $|M|\ne0$.

Let $\mathbf{Q}_0=\mathbf{p}_2$ and $\mathbf{Q}_1=\mathbf{p}_1+\mathbf{p}_2t$. Taking the cross product of $\mathbf{Q}_i$ and $\mathbf{F}(t)$ gives $\mathbf{Q}_i\times\mathbf{F}(t)$, $i=0,1$.
Considering the cross product of $\mathbf{Q}_1\times\mathbf{F}(t)$ and $\mathbf{Q}_0\times\mathbf{F}(t)$, we have
\[
\bigl(\mathbf{Q}_1\times\mathbf{F}(t)\bigr)\times\bigl(\mathbf{Q}_0\times\mathbf{F}(t)\bigr)=|\mathbf{p}_0\mathbf{p}_1\mathbf{p}_2|\mathbf{F}(t),
\]
and $\mathrm{deg}\left(\mathbf{Q}_1\times\mathbf{F}(t)\right)+\mathrm{deg}\left(\mathbf{Q}_0\times\mathbf{F}(t)\right)=\mathrm{deg}(\mathbf{F}(t))$.
Therefore, $\mathbf{Q}_1\times\mathbf{F}(t)$ and $\mathbf{Q}_0\times\mathbf{F}(t)$ form a $\mu$-basis of the curve $\mathbf{F}(t)$, which is consistent with the result of Wang \cite{wang2016explicit_mu_bases}.

\subsection{Cubic Non-degenerate Curves}\label{subsec:4.2}

Let $\mathbf{F}(t)=\mathbf{p}_0+\mathbf{p}_1t+\mathbf{p}_2t^2+\mathbf{p}_3t^3$ be a curve, where $\mathbf{p}_i\in\mathbb{R}^3,i=0,1,2,3$, $\mathbf{p}_0\ne0$ and $\mathbf{p}_3\ne0$.
Let $M=(\mathbf{p}_0^T,\mathbf{p}_1^T,\mathbf{p}_2^T,\mathbf{p}_3^T)^T$ be a $4\times3$ matrix, and $Rank(M)=3$.

Let $\mathbf{Q}_0=\mathbf{p}_3$, $\mathbf{Q}_1=\mathbf{p}_2+\mathbf{p}_3t$, and $\mathbf{Q}_2=\mathbf{p}_1+\mathbf{p}_2t+\mathbf{p}_3t^2$. Taking the cross product of $\mathbf{Q}_i$ and $\mathbf{F}(t)$ gives
$\mathbf{Q}_i\times\mathbf{F}(t)$, $i=0,1,2$. Considering the cross product of $\mathbf{Q}_2\times\mathbf{F}(t)$ and $\mathbf{Q}_0\times\mathbf{F}(t)$, we have
\begin{eqnarray}\label{eq:4.1}
	\bigl(\mathbf{Q}_2\times\mathbf{F}(t)\bigr)\times\bigl(\mathbf{Q}_0\times\mathbf{F}(t)\bigr)=\left[\frac{(\mathbf{p}_3\times\mathbf{p}_0)\cdot\mathbf{F}(t)}{t}\right]\mathbf{F}(t)=|\mathbf{p}_0\mathbf{p}_1\mathbf{p}_3|+|\mathbf{p}_0\mathbf{p}_2\mathbf{p}_3|t.
\end{eqnarray}

1)When $\mathbf{p}_3\times\mathbf{p}_0=\mathbf{0}$，$\mathbf{Q}_2\times\mathbf{F}(t)$ and $\mathbf{Q}_0\times\mathbf{F}(t)$ are linearly dependent over $\mathbb{R}[t]$.
In this case, considering the cross product of $\mathbf{Q}_1\times\mathbf{F}(t)$ and $\mathbf{Q}_0\times\mathbf{F}(t)$, we have
\[
\bigl(\mathbf{Q}_1\times\mathbf{F}(t)\bigr)\times\bigl(\mathbf{Q}_0\times\mathbf{F}(t)\bigr)=|\mathbf{p}_1\mathbf{p}_2\mathbf{p}_3|\mathbf{F}(t),
\]
and $\mathrm{deg}\left(\mathbf{Q}_1\times\mathbf{F}(t)\right)+\mathrm{deg}\left(\mathbf{Q}_0\times\mathbf{F}(t)\right)=\mathrm{deg}(\mathbf{F}(t))$，
So $\mathbf{Q}_1\times\mathbf{F}(t)$ and $\mathbf{Q}_0\times\mathbf{F}(t)$ form a $\mu$-basis of the curve $\mathbf{F}(t)$.

2)When $\mathbf{p}_3\times\mathbf{p}_0\ne\mathbf{0}$，$\eqref{eq:4.1}$ holds, and $\mathbf{Q}_2\times\mathbf{F}(t)$ and $\mathbf{Q}_0\times\mathbf{F}(t)$
are linearly independent over $\mathbb{R}[t]$. If $\mathbf{Q}_2\times\mathbf{F}(t)$ can be decomposed as
\[
\mathbf{Q}_2\times\mathbf{F}(t)=\left(|\mathbf{p}_0\mathbf{p}_1\mathbf{p}_3|+|\mathbf{p}_0\mathbf{p}_2\mathbf{p}_3|t\right)\cdot\left(\frac{\mathbf{p}_1\times\mathbf{p}_0}{|\mathbf{p}_0\mathbf{p}_1\mathbf{p}_3|}+\frac{\mathbf{p}_3\times\mathbf{p}_0}{|\mathbf{p}_0\mathbf{p}_2\mathbf{p}_3|}t\right),
\]
it is required that $\displaystyle\frac{|\mathbf{p}_0\mathbf{p}_1\mathbf{p}_3|^2}{|\mathbf{p}_0\mathbf{p}_2\mathbf{p}_3|}=|\mathbf{p}_0\mathbf{p}_1\mathbf{p}_2|$.
In this case, let $\mathbf{L}(t)=\displaystyle\frac{\mathbf{p}_1\times\mathbf{p}_0}{|\mathbf{p}_0\mathbf{p}_1\mathbf{p}_3|}+\frac{\mathbf{p}_3\times\mathbf{p}_0}{|\mathbf{p}_0\mathbf{p}_2\mathbf{p}_3|}t$,
then we have
\[
\mathbf{L}(t)\times\bigl(\mathbf{Q}_0\times\mathbf{F}(t)\bigr)=\mathbf{F}(t)\,\,\,,\,\,\mathrm{deg}\left(\mathbf{L}(t)\right)+\mathrm{deg}\left(\mathbf{Q}_0\times\mathbf{F}(t)\right)=\mathrm{deg}(\mathbf{F}(t)),
\]
So $\mathbf{L}(t)$ and $\mathbf{Q}_0\times\mathbf{F}(t)$ form a $\mu$-basis of the curve $\mathbf{F}(t)$.

When $\displaystyle\frac{|\mathbf{p}_0\mathbf{p}_1\mathbf{p}_3|^2}{|\mathbf{p}_0\mathbf{p}_2\mathbf{p}_3|}\ne|\mathbf{p}_0\mathbf{p}_1\mathbf{p}_2|$，according to 
$$
|\mathbf{p}_1\mathbf{p}_2\mathbf{p}_3|(\mathbf{p}_3\times \mathbf{p}_0)-|\mathbf{p}_0\mathbf{p}_2\mathbf{p}_3|(\mathbf{p}_3\times \mathbf{p}_1)+|\mathbf{p}_0\mathbf{p}_1\mathbf{p}_3|(\mathbf{p}_3\times \mathbf{p}_2)=0,
$$
we construct $\mathbf{L}(t)=|\mathbf{p}_1\mathbf{p}_2\mathbf{p}_3|(\mathbf{Q}_2\times \mathbf{F}(t))-|\mathbf{p}_0\mathbf{p}_2\mathbf{p}_3|(\mathbf{Q}_1\times \mathbf{F}(t))+|\mathbf{p}_0\mathbf{p}_1\mathbf{p}_3|(\mathbf{Q}_0\times \mathbf{F}(t))$,Expanding it gives $\mathbf{L}(t)=\mathbf{L}_0+\mathbf{L}_1 t$,where 
\begin{eqnarray*}
	\mathbf{L}_0&=&|\mathbf{p}_1\mathbf{p}_2\mathbf{p}_3|(\mathbf{p}_1\times \mathbf{p}_0)-|\mathbf{p}_0\mathbf{p}_2\mathbf{p}_3|(\mathbf{p}_2\times \mathbf{p}_0)+|\mathbf{p}_0\mathbf{p}_1\mathbf{p}_3|(\mathbf{p}_3\times \mathbf{p}_0),\\
	\mathbf{L}_1&=&-|\mathbf{p}_0\mathbf{p}_2\mathbf{p}_3|(\mathbf{p}_3\times \mathbf{p}_0)+|\mathbf{p}_0\mathbf{p}_1\mathbf{p}_3|(\mathbf{p}_3\times \mathbf{p}_1)-|\mathbf{p}_0\mathbf{p}_1\mathbf{p}_2|(\mathbf{p}_3\times \mathbf{p}_2).
\end{eqnarray*}
It satisfies
\[
\begin{cases}
	\mathbf{L}(t)\cdot\mathbf{F}(t)=0,\\
	\left(\mathbf{Q}_2\times\mathbf{F}(t)\right)\times\mathbf{L}(t)=\left(|\mathbf{p}_0\mathbf{p}_1\mathbf{p}_3|^2-|\mathbf{p}_0\mathbf{p}_1\mathbf{p}_2|\cdot|\mathbf{p}_0\mathbf{p}_2\mathbf{p}_3|\right)\cdot \mathbf{F}(t).
\end{cases}
\]
Therefore
\begin{eqnarray*}
	\mathbf{L}(t)&=&-|\mathbf{p}_1\mathbf{p}_2\mathbf{p}_3|(\mathbf{p}_1\times \mathbf{p}_0)+|\mathbf{p}_0\mathbf{p}_2\mathbf{p}_3|(\mathbf{p}_2\times \mathbf{p}_0)-|\mathbf{p}_0\mathbf{p}_1\mathbf{p}_3|(\mathbf{p}_3\times \mathbf{p}_0)\\
	&&+\left[|\mathbf{p}_0\mathbf{p}_2\mathbf{p}_3|(\mathbf{p}_3\times \mathbf{p}_0)-|\mathbf{p}_0\mathbf{p}_1\mathbf{p}_3|(\mathbf{p}_3\times \mathbf{p}_1)+|\mathbf{p}_0\mathbf{p}_1\mathbf{p}_2|(\mathbf{p}_3\times \mathbf{p}_2)\right]t,\\
	\mathbf{Q}_2\times\mathbf{F}(t)&=&(\mathbf{p}_1\times\mathbf{p}_0)+(\mathbf{p}_2\times \mathbf{p}_0)t+(\mathbf{p}_3\times\mathbf{p}_0)t^2.
\end{eqnarray*}
form a $\mu$-basis of the curve $\mathbf{F}(t)$, which is consistent with the result of Wang \cite{wang2016explicit_mu_bases}.

\section{Further discussion}\label{sec:5}
In this paper, by constructing the vector polynomials $\mathbf{Q}_i\times\mathbf{F}(t)$, 
we give explicit formulas for the $\mu$-bases of quartic non-degenerate real plane rational curves, 
and use these $\mu$-bases to derive closed-form expressions for their implicit equations. 
In the future, how to derive explicit formulas for the $\mu$-bases of plane rational curves of arbitrary degree 
will be a problem that we need to further study and solve.

\end{document}